\documentclass[12pt]{article}
\usepackage{amsmath,amssymb,amstext,dsfont,fancyvrb,float,fontenc,graphicx,subfigure,theorem,hyperref}
\usepackage[utf8]{inputenc}

\usepackage{tikz,everypage-1x}

\usepackage[letterpaper]{geometry}
\ifx\volno\undefined\def\volno{0}\fi
\ifx\volyear\undefined\def\volyear{2017}\fi
\ifx\pagno\undefined\def\pagno{000--000}\fi

\newfont{\footsc}{cmcsc10 at 8truept}
\newfont{\footbf}{cmbx10 at 8truept}
\newfont{\footrm}{cmr10 at 10truept}

\usepackage{fancyhdr}
\usepackage{relsize}
\usepackage{sectsty}
\allsectionsfont{\larger[-1]} 

\renewcommand\paragraph{\@startsection{paragraph}{4}{\z@}
                                    {2ex \@plus.5ex \@minus.2ex}
                                    {-1em}
                                    {\normalfont\normalsize\bfseries}}

\renewcommand\subparagraph{\@startsection{subparagraph}{5}{\parindent}
                                       {2ex \@plus.5ex \@minus .2ex}
                                       {-1em}
                                      {\normalfont\normalsize\bfseries}}

\newlength{\BiblioSpacing}
\renewenvironment{thebibliography}[1]{
\begin{oldthebibliography}{#1}
\setlength{\parskip}{\BiblioSpacing}
\setlength{\itemsep}{\BiblioSpacing}
}
{
\end{oldthebibliography}
}

\usepackage[strict]{changepage}
\def\abstractname{Abstract -}   
\def\abstract{\begin{adjustwidth}{1cm}{1cm} \par    \footnotesize \noindent {\bf \abstractname} 
\def\endabstract{ \end{adjustwidth} \smallskip }}

{\theorembodyfont{\itshape}\newtheorem{theorem}{Theorem}[section]}
{\theorembodyfont{\itshape}\newtheorem{proposition}[theorem]{Proposition}}
{\theorembodyfont{\itshape}\newtheorem{definition}[theorem]{Definition}}
{\theorembodyfont{\itshape}\newtheorem{lemma}[theorem]{Lemma}}
{\theorembodyfont{\itshape}\newtheorem{corollary}[theorem]{Corollary}}
{\theorembodyfont{\rm}}
{\theorembodyfont{\rm}\newtheorem{remark}[theorem]{Remark}}
{\theorembodyfont{\rm}}
{\theorembodyfont{\rm }\newtheorem{example}[theorem]{Example}}

\usepackage[all]{xy}
\theoremstyle{plain}
\numberwithin{equation}{section}
\numberwithin{figure}{section}
\numberwithin{table}{section}

\newenvironment{solution}{\textit{Proof.}}{\hfill$\blacksquare$}

\title{\Large\bf Products of $p$-Adic Valuation Trees}
\author{\sc D. Snyder\thanks{This work was supported by an REU research grant (NSF grant No. 1851948)}}

\begin{document}
\setcounter{page}{1}
\maketitle
\thispagestyle{fancy}

\vskip 1.5em

\begin{abstract}
The study of prime divisibility plays a crucial role in number theory. The $p$-adic valuation of a number is the highest power of a prime, $p$, that divides that number. Using this valuation, we construct $p$-adic valuation trees to visually represent the valuations of a sequence. We investigate how nodes split on trees generated by linear functions with rational coefficients, as well as those formed from a product of linear or lower degree polynomials. We describe the infinite branches of these polynomial trees and the valuations of their terminating nodes.
\end{abstract}
 

\begin{MSC}
11B83; 11S05; 11Y55
\end{MSC}

\section{Introduction}
A core aspect of number theory is investigating patterns within primes and their divisibility. One way we measure this divisibility is with the $p$-adic valuation, the greatest power of the prime $p$ which divides a given integer. This definition is extended to rational numbers and provides a new way to categorize the field of rationals. Using a sequence of rational numbers, we can create a $p$-adic valuation sequence which describes the $p$-adic valuation of each term. These sequences show a connection between a function and how the powers of a prime number appear in its outputs, as Guan \cite{guan} notes in Chapter 5 of their thesis. To better recognize the patterns within these $p$-adic valuation sequences, we construct $p$-adic valuation trees that organize the terms into nodes and create a visual aid from which to view these sequences. The structure of $p$-adic valuation trees and their relation to their generating sequences could provide insight into features of prime divisibility. For this paper, we will discuss sequences generated by polynomial functions.

Previously, the works from Almodovar et al. \cite{almo}, Boultinghouse et al. \cite{2adic}, Kozhushkina et al. \cite{n2+7}, and Nguyen \cite{nguyen} characterize $2$-adic valuation sequences and trees for specific quadratic polynomials, depending upon their coefficients. Some polynomials discussed are $n^2+a$ for some integer, $a$, or more general quadratics with integer coefficients, written as $an^2+bn+c$. However, these characterizations are not complete.  Boultinghouse et al. \cite{padic} also classified the structure of $p$-adic valuation trees for all quadratic polynomials using their discriminant for odd primes $p$. Beyond quadratic polynomials, Guan \cite{guan} conjectured about higher-degree polynomials and their $p$-adic valuation trees.

In this paper, we will describe the exact structure of a $p$-adic valuation tree generated by a linear polynomial with rational coefficients. Then, we will characterize the $p$-adic valuation trees for quadratic polynomials which are factorable over the rationals. Next, we construct, in its entirety, the $p$-adic valuation trees for any completely factorable polynomial with rational coefficients. Finally, we remark upon polynomials which factor into nonlinear parts and their $p$-adic valuation trees.

\section{Background}
\subsection{p-Adic numbers}
\begin{definition}
    Given a prime $p$, the $p$-adic valuation of an integer $n$, denoted $\nu_p(n)$, is an integer $x$, such that \[p^x\mid n \quad \text{and} \quad p^{x+1}\nmid n.\] This can be extended to find the $p$-adic valuation of any rational number $r$. Written in the form $r=\frac{a}{b}$, the valuation satisfies the property \[\nu_p(r)=\nu_p\left(\frac{a}{b}\right)=\nu_p(a)-\nu_p(b).\]
\end{definition}


\begin{remark}\label{val.prod}
    If $n=ab$ for $a,b\in \mathbb{Q}$, then $\nu_p(n)=\nu_p(a)+\nu_p(b)$.
\end{remark}

Using the $p$-adic valuation, we create an absolute value on the rational numbers.

\begin{definition}
    We define the $p$-adic absolute value, $|\cdot |_p: \mathbb{Q}\to \mathbb{Q}$, to be \[|x|_p=p^{-\nu_p(x)}\] for any $x\in \mathbb{Q}-\{0\}$. If $x=0$, we say $|0|_p=0$.
\end{definition}


With the $p$-adic absolute value, we now define the field of $p$-adic numbers which will be crucial to this paper. We use \cite{Gouvea} as a reference for the following definition and theorem.

\begin{definition}
    The $p$-adic numbers are a completion of the rational numbers with respect to the $p$-adic absolute value. This is denoted by $\mathbb{Q}_p$.
\end{definition}

\begin{theorem}
Every $p$-adic number $\alpha \in \mathbb{Q}_p$, can be uniquely written as the series
\begin{equation*}
 \alpha = \sum^\infty_{n=n_0} a_np^n =a_{n_0}p^{n_0}+a_{n_0+1}p^{n_0+1}+ \cdots,    
\end{equation*}
where $n_0\in \mathbb{Z} \cup \{ \infty \}, a_n\in \mathbb{Z}, 0 \leq a_n \leq p-1, \text{ and }
a_{n_0} \neq 0.$ When $n_0\geq 0$, the number $\alpha$ is a $p$-adic integer. The set of $p$-adic integers is denoted by $\mathbb{Z}_p$.
\end{theorem}


Notice that the $p$-adic integers can be described using modular arithmetic. Consider an arbitrary $p$-adic integer \[\alpha=a_{0}+a_{1}p+ \cdots+ a_np^n +\cdots \] where $0\leq a_n\leq p-1$. Then we have \begin{align*}
    \alpha&\equiv a_0 \bmod{p}\\
    &\equiv a_0+a_1p \bmod{p^2}\\
    & \hspace{2mm} \vdots \\
    &\equiv a_{0}+a_{1}p+ \cdots+ a_np^n \bmod{p^{n+1}}.
\end{align*} We can also construct $p$-adic valuation sequences by finding the $p$-adic valuation of each term in a given sequence. Throughout this paper, our sequences will be generated by polynomial functions.

\begin{example}\label{n^2+4}
    Consider the sequence generated by $f(n)=n^2+4$: \[\{f(n)\}_{n\geq 0}=\{4, 5, 8, 13, 20, 29, 40, \dots\}.\] The corresponding $2$-adic valuation sequence is \[\{\nu_2(f(n))\}_{n\geq 0}=\{2, 0, 3, 0, 2, 0, 3, \dots\}.\]
\end{example}

Next, we will use these $p$-adic valuation sequences to construct valuation trees.







\subsection{p-Adic valuation trees}
Consider a sequence of rational numbers, $\{S_n\}_{n\geq 0}$, and its $p$-adic valuations, $\{\nu_p(S_n)\}_{n\geq 0}$. We can construct a tree to visually represent this valuation sequence, using nodes and branches. The node on the top of the tree contains all the terms in the sequence. If this node splits, it travels to the next level with $p$ branches, each containing the terms in the sequence of the form $c_0+np$ where $0\leq c_0\leq p-1$. For this paper, the branches will be organized from left to right with $c_0$ increasing from $0$ to $p-1$. These branches end at the node on the $1^{st}$ level. At the $k^{th}$ level, if a node splits, its branches will contain the terms in the sequence of the form $c_0+c_1p+\cdots +c_{k-1}p^{k-1}+np^k$, where $c_{i}$ has been previously determined for all $0\leq i\leq k-1$. In the tree, the number $c_0+c_1p+\cdots +c_{k-1}p^{k-1}$ will be written above its respective branch.

Now let us define the conditions by which a node will split.
\begin{definition}
    A node in the $p$-adic valuation tree of $\{\nu_p(S_n)\}_{n\geq 0}$ is called \textit{terminating} when the valuation is the same for all terms in that node. In the tree, this will be denoted by a circle with the valuation written inside. This type of node will not split.
\end{definition}

\begin{definition}
    A node is called \textit{non-terminating} when at least two terms within the same node have different valuations. In the tree, this will be denoted by a circle with an asterisk inside. This node will split into the next level.
\end{definition}

The sequences discussed in this paper will be constructed from polynomials with integer or rational coefficients, however, these trees can be constructed from other sequences. Following the work done by Boultinghouse et al. \cite{padic} and others, we will focus on polynomials.

\begin{example}
    The figure below is the $2$-adic valuation tree representing $\{\nu_p(S_n)\}_{n\geq 0}$, where $S_n=n^2+4$: 
\[\xymatrix{
 & & *+[Fo]{*} \ar@{-}[dl]_{0} \ar@{-}[dr]^{1} & \\
 & *+[Fo]{*} \ar@{-}[dl]_{0} \ar@{-}[dr]^{2} & & *+[Fo]{0} \\
 *+[Fo]{2} & & *+[Fo]{3} &
 } \]
 \[\text{Figure 2.1: $2$-Adic valuation tree of $n^2+4$}\]
This is the visual representation of the $2$-adic valuation sequence in Example \ref{n^2+4}. This tree has $3$ levels. We label the highest level as the $0^{th}$ level. The $0^{th}$ level has a non-terminating node. This node splits into two branches, one where the terms have the form $2n$ and another with the form $1+2n$. Hence, $0$ is written above the branch to the left and $1$ is written above the branch to the right. The $1^{st}$ level has a non-terminating node of the form $2n$ and a terminating node of the form $1+2n$ with a valuation of $0$. The non-terminating node of the form $2n$ splits into two branches, one where the terms have the form $2(2n)=4n$ and another with the form $2(1+2n)=2+4n$. Hence, $0$ is written above the branch to the left and $2$ is written above the branch to the right. The $2^{nd}$ level has a terminating node of the form $4n$ with a valuation of $2$ and a terminating node of the form $2+4n$ with a valuation of $3$. This also shows us that the non-terminating node at the $1^{st}$ level has a valuation of at least $2$ due to the fixed valuations of the nodes beneath it.
\end{example}

\begin{definition}
    A tree is \textit{infinite} if every level has a non-terminating node. A tree is \textit{finite} if it is not infinite.
\end{definition}

The previous example is a finite tree because it does not have a non-terminating node at the $2^{nd}$ level.

Since every non-terminating node is contained within a non-terminating node on the previous level, an infinite tree will have at least one path of non-terminating nodes that passes through every level. We will refer to this path as an \textit{infinite branch}. For an infinite branch, every level will have a non-terminating node which will split. In addition, one of the successive nodes will also be non-terminating. This creates an infinite path from the $0^{th}$ level down through every other level using these non-terminating nodes. The following theorem shows a useful connection between infinite trees, infinite branches, and their generating polynomial.

\begin{theorem}(L. Almodovar et al. \cite{almo})\label{almo}
    Any infinite branch in the tree associated to the polynomial $f$ corresponds to a root of $f(x)=0$ in the $p$-adic field $\mathbb{Q}_p$.
\end{theorem}

\begin{corollary}(L. Almodovar et al. \cite{almo})
    The $p$-adic valuation $\nu_p(f(n))$ admits a closed-form formula if
the equation $f(x) = 0$ has no solutions in $\mathbb{Q}_p$.
\end{corollary}

However, if a root of a polynomial is a $p$-adic number, it does not imply that this corresponds to an infinite branch on the tree for this polynomial, nor does this tree have to be infinite. An example of this is shown in Remark \ref{zero point}. The root of a polynomial will only correspond to an infinite branch on the tree for this polynomial if that root is a $p$-adic integer. More specifically, if a polynomial has a root in the $p$-adic integers, then its respective tree will be infinite. The reasoning behind this will be seen throughout the proofs of this paper. We will explore this in the next section in a simple case.

\section{Linear trees}

First, we will describe the $p$-adic valuation tree for a linear polynomial with integer coefficients.

\begin{proposition}\label{lin.descr}
    Let $f(n)=an+b\in \mathbb{Z}[n]$, where $a\neq 0$ and $\gcd(a,b)=1$. If $-\frac{b}{a}\in \mathbb{Z}_p$, then a non-terminating node at the $k^{th}$ level splits into $p-1$ terminating nodes with a valuation of $k$ and one non-terminating node with a valuation of at least $k+1$.
\end{proposition}

\begin{solution}
    We will prove this by induction. Write the $p$-adic expansion for $-\frac{b}{a}$ as \[-\frac{b}{a}=c_0+c_1p+c_2p^2+\cdots,\] where $0\leq c_i \leq p-1$. First, consider the $0^{th}$ level, where the non-terminating node splits into $p$ branches of the form $j_0+np$ for $0\leq j_0\leq p-1$. Notice that the congruence \begin{align*}
        f(j_0+np)&= a(j_0+np)+b\\
        &\equiv aj_0+b \bmod{p}\\
        &\equiv 0 \bmod{p}
    \end{align*} has the solution $j_0\equiv -\frac{b}{a}\equiv c_0 \bmod{p}.$ Thus, when $j_0\not\equiv c_0 \bmod{p}$, the branch has a terminating node with a valuation of $0$. When $j_0\equiv c_0 \bmod{p}$, the branch has a valuation of at least $1$. To see that this branch has a non-terminating node, consider \[c_0+c_1p+c_2p^2+\cdots +c_sp^s\] with $s$ large enough. This natural number is within the node and is arbitrarily close to the root $-\frac{b}{a}$. Since $f(n)$ is a continuous function, $f(c_0+c_1p+c_2p^2+\cdots +c_sp^s)$ is arbitrarily close to $f\left(-\frac{b}{a}\right)=0$. Thus, \[\nu_p(f(c_0+c_1p+c_2p^2+\cdots +c_sp^s))\] is arbitrarily large. Now, consider the natural number \[c_0+(c_1+1)p+c_2p^2+\cdots +c_sp^s,\] which is also within the same node but not as close to $-\frac{b}{a}$. Again, by the continuity and linearity of $f(n)$, $f(c_0+(c_1+1)p+c_2p^2+\cdots +c_sp^s)$ is not as close to $f\left(-\frac{b}{a}\right)=0$, so \[\nu_p(f(c_0+(c_1+1)p+c_2p^2+\cdots +c_sp^s))\] is some fixed finite value that is different than the previous valuation. Hence, the branch has two elements with different valuations, signifying a non-terminating node.

    Next, suppose that the statement is true for all levels up to the $k^{th}$ level. For the non-terminating node at the $k^{th}$ level, branches to the $(k+1)^{st}$ level have the form \[N_{k+1}=c_0+c_1p+\cdots+c_{k-1}p^{k-1}+j_kp^k+np^{k+1},\] where $0\leq j_k\leq p-1$. By the inductive step, we know that \[\nu_p(f(N_{k+1}))\geq k,\] and so $f(N_{k+1})\equiv 0 \bmod{p^k}$. Notice that the congruence \begin{align*}
        f(N_{k+1})&= a(N_{k+1})+b\\
        &\equiv a(c_0+c_1p+\cdots+c_{k-1}p^{k-1}+j_kp^k+np^{k+1})+b \bmod{p^{k+1}}\\
        &\equiv 0 \bmod{p^{k+1}}
    \end{align*} has the solution \[c_0+c_1p+\cdots+c_{k-1}p^{k-1}+j_kp^k\equiv -\frac{b}{a}\equiv c_0+c_1p+\cdots+c_{k-1}p^{k-1}+c_kp^k \bmod{p^{k+1}},\] which implies $j_k\equiv c_k \bmod{p}$. When $j_k\not\equiv c_k \bmod{p}$, then $f(N_{k+1})\not\equiv 0 \bmod{p^{k+1}}$ and so $\nu_p(f(N_{k+1}))=k$. Thus, the branch has a terminating node with a valuation of $k$. When $j_k\equiv c_k \bmod{p}$, then $f(N_{k+1})\equiv 0 \bmod{p^{k+1}}$ and so $\nu_p(f(N_{k+1}))\geq k+1$. Hence, this node has a valuation of at least $k+1$. To see that this node is non-terminating, consider \[c_0+c_1p+c_2p^2+\cdots c_{k+1}p^{k+1}+\cdots+c_sp^s\] with $s$ large enough. This natural number is within the node and is arbitrarily close to the root $-\frac{b}{a}$. Since $f(n)$ is a continuous function, $f(c_0+c_1p+c_2p^2+\cdots c_{k+1}p^{k+1}+\cdots+c_sp^s)$ is arbitrarily close to $f\left(-\frac{b}{a}\right)=0$. Thus, \[\nu_p(f(c_0+c_1p+c_2p^2+\cdots +c_{k+1}p^{k+1}+\cdots +c_sp^s))\] is arbitrarily large. Now, consider the natural number \[c_0+c_1p+c_2p^2+\cdots+ (c_{k+1}+1)p^{k+1}+\cdots +c_sp^s,\] which is also within the same node but not as close to $-\frac{b}{a}$. Again, by the continuity and linearity of $f(n)$, $f(c_0+c_1p+c_2p^2+\cdots+ (c_{k+1}+1)p^{k+1}+\cdots +c_sp^s)$ is not as close to $f\left(-\frac{b}{c}\right)=0$ so \[\nu_p(f(c_0+c_1p+c_2p^2+\cdots+ (c_{k+1}+1)p^{k+1}+\cdots +c_sp^s))\] is some fixed finite value that is different than the previous valuation. Hence, the branch has two elements with different valuations, signifying a non-terminating node.

    Thus, through induction, we see that for any level $k$, there is a non-terminating node that splits into $p-1$ terminating nodes with a valuation of $k$ and a non-terminating node with a valuation of at least $k+1$.
\end{solution}

\vspace{3mm}

Thus, the $p$-adic valuation tree for a linear function with a $p$-adic integer root has one infinite branch. Additionally, notice that the infinite branch corresponds to the $p$-adic expansion for $-\frac{b}{a}$ as the proof also shows that $j_r\equiv c_r \bmod{p^r}$ for all $r\geq 0$.

\begin{example}\label{2n-5}
Now, consider when $p=3$ and $f(n)=2n-5$. We can write the root of this polynomial as the $p$-adic expansion \[\frac{5}{2}= 1\cdot 3^0+2\cdot 3^1+1\cdot 3^2+1\cdot 3^3+1\cdot 3^4+\cdots\] and see that this is a $3$-adic integer. The tree for $\{\nu_3(2n-5)\}_{n\geq 0}$ is shown below. 
\[
\xymatrix{
 & & *+[Fo]{*} \ar@{-}[dl]_{0} \ar@{-}[d]^{1} \ar@{-}[dr]^{2} & & \\
 &  *+[Fo]{0} & *+[Fo]{*} \ar@{-}[dl]_{1} \ar@{-}[d]^{4} \ar@{-}[dr]^{7} & *+[Fo]{0} & \\
 &  *+[Fo]{1} & *+[Fo]{1} & *+[Fo]{*} \ar@{-}[dl]_{7} \ar@{-}[d]^{16} \ar@{-}[dr]^{25} & \\
 & & *+[Fo]{2} & *+[Fo]{*} \ar@{-}[dl]_{16} \ar@{-}[d]^{43} \ar@{-}[dr]^{70} & *+[Fo]{2} \\
 & & *+[Fo]{3} & *+[Fo]{*} \ar@{-}[dl]_{43} \ar@{-}[d]^{124} \ar@{-}[dr]^{205} & *+[Fo]{3} \\ 
 & &  *+[Fo]{4} & *+[Fo]{*} & *+[Fo]{4}
}
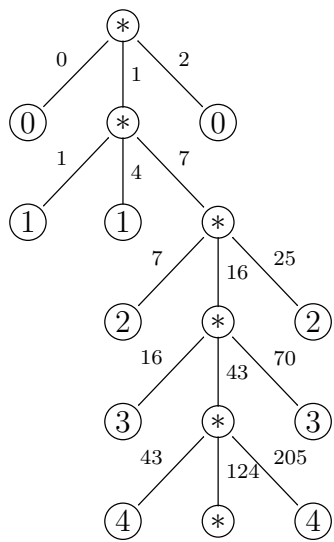
\] \[\text{Figure 3.1: $3$-Adic valuation tree for $2n-5$}\] The infinite branch matches the $3$-adic expansion for $\frac{5}{2}$, following our hypothesis.
\end{example}

\begin{remark}\label{zero point}
    For the case where $-\frac{b}{a}\notin \mathbb{Z}_p$, it must be true that $a\equiv 0 \bmod{p}$ and $b\not\equiv 0\bmod{p}$. We can write $a=pr$ and $b=ps+t$ for $r,s,t\in \mathbb{Z}$, where $1\leq t\leq p-1$. Thus, every term in the sequence $\{f(n)\}_{n\geq 0}$ has the form \begin{align*}
        an+b&=prn+ps+t\\
        &=p(rn+s)+t\\
        &\equiv t \bmod{p}\\
        &\not\equiv 0 \bmod{p},
    \end{align*} and so $\nu_p(f(n))=0$. Hence, the tree consists of a single terminating node with a valuation of $0$.
\end{remark}

\section{Combining two linear trees}
First, let us discuss the effect of combining a linear tree with a tree generated by a constant.

\begin{lemma}\label{coef}
    Let $f(n)=c(an+b)$, where $a,b\in \mathbb{Z}$, $\gcd(a,b)=1$, $c\in \mathbb{Q}$, and $a,c\neq 0$. Consider the tree generated by $\{\nu_p(f(n))\}_{n\geq 0}$. If $-\frac{b}{a}\in \mathbb{Z}_p$, then a non-terminating node at the $k^{th}$ level splits into $p-1$ terminating nodes with a valuation of $\nu_p(c)+k$ and one non-terminating node with a valuation of at least $\nu_p(c)+k+1$. If $-\frac{b}{a}\notin \mathbb{Z}_p$, then the tree is a singular terminating node with a valuation of $\nu_p(c)$.
\end{lemma}

\begin{solution}
    The proof is trivial given the property of $p$-adic valuations seen in Remark \ref{val.prod}. This implies that $\nu_p(f(n))=\nu_p(c)+\nu_p(an+b)$.
\end{solution}

\vspace{3mm}

Notice that because $c\in \mathbb{Q}$, we have described the trees for every linear polynomial with rational coefficients.


\begin{example}\label{3/2n+9}
 Consider when $p=3$ and $f(n)=\frac{3}{2}(n+6)=\frac{3}{2}n+9$. We can write the root of this polynomial as the $p$-adic expansion \[-6=0\cdot 3^0 +1\cdot 3^1+2\cdot 3^2+2\cdot 3^3+2\cdot 3^4 + \cdots\] and see that this is a $3$-adic integer. Since $\nu_3\left(\frac{3}{2}\right)=1$, the tree for $\{\nu_3\left(\frac{3}{2}n+9\right)\}_{n\geq 0}$ is as follows.
\[
\xymatrix{
 & & *+[Fo]{*} \ar@{-}[dl]_{0} \ar@{-}[d]^{1} \ar@{-}[dr]^{2} & & \\
 & *+[Fo]{*}  \ar@{-}[dl]_{0} \ar@{-}[d]^{3} \ar@{-}[dr]^{6} & *+[Fo]{1} & *+[Fo]{1} & \\
 *+[Fo]{2} & *+[Fo]{*} \ar@{-}[dl]_{3} \ar@{-}[d]^{12} \ar@{-}[dr]^{21} & *+[Fo]{2} & & \\
 *+[Fo]{3} & *+[Fo]{3} & *+[Fo]{*} \ar@{-}[dl]_{21} \ar@{-}[d]^{48} \ar@{-}[dr]^{75} & & \\
 & *+[Fo]{4} & *+[Fo]{4} & *+[Fo]{*}  \ar@{-}[dl]_{75} \ar@{-}[d]^{156} \ar@{-}[dr]^{237} & \\
 & & *+[Fo]{5} & *+[Fo]{5} & *+[Fo]{*}
}
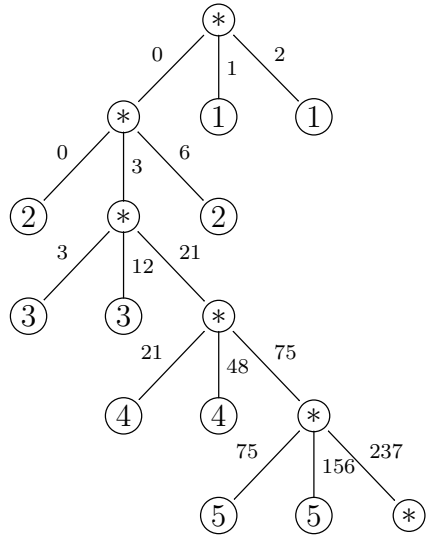
\] \[\text{Figure 4.1: $3$-Adic valuation tree for $\frac{3}{2}n+9$}\] Thus, we see that the infinite branch matches the $3$-adic expansion $-6$. The valuations of the terminating nodes are one greater due to the factor of $\frac{3}{2}$.
\end{example}

The following theorem will describe the $p$-adic valuation tree for the product of two linear polynomials in $\mathbb{Z}[n]$.

\begin{theorem}
    Let $f(n)=(a_1n+b_1)(a_2n+b_2)$, where $a_i,b_i\in \mathbb{Z}$, $a_i\neq 0$, and $\gcd(a_i,b_i)=1$ for $i\in \{1,2\}$. Consider the roots of the polynomial $-\frac{b_1}{a_1}$ and $-\frac{b_2}{a_2}$. We have the following cases: \begin{itemize}
        \item[(i)] If $-\frac{b_1}{a_1}= -\frac{b_2}{a_2}$, then the tree generated by $\{\nu_p(f(n))\}_{n\geq 0}$ has the same terminating and non-terminating nodes as the tree generated by $\{\nu_p(a_1n+b_1)\}_{n\geq 0}$. However, the valuations at these nodes are twice what they are for $\{\nu_p(a_1n+b_1)\}_{n\geq 0}$.
        \item[(ii)] If $-\frac{b_2}{a_2}\notin \mathbb{Z}_p$, then the tree generated by $\{\nu_p(f(n))\}_{n\geq 0}$ has the same terminating and non-terminating nodes as the tree generated by $\{\nu_p(a_1n+b_1)\}_{n\geq 0}$. Moreover, the terminating nodes have the same valuations as those on the tree generated by $\{\nu_p(a_1n+b_1)\}_{n\geq 0}$.
        \item[(iii)] If $-\frac{b_1}{a_1},-\frac{b_2}{a_2}\in \mathbb{Z}_p$ and $-\frac{b_1}{a_1}\neq -\frac{b_2}{a_2}$, then either there exists an integer $r$ such that $-\frac{b_1}{a_1}\equiv -\frac{b_2}{a_2} \bmod{p^r}$ and $-\frac{b_1}{a_1}\not\equiv -\frac{b_2}{a_2} \bmod{p^{r+1}},$ or we say $r=0$. When $0\leq k\leq \max{\{r-1, 0\}}$, the non-terminating node at the $k^{th}$ level splits into $p-1$ terminating nodes with a valuation of $2k$ and one non-terminating node with a valuation of at least $2(k+1)$. At the $r^{th}$ level, the non-terminating node splits into $p-2$ terminating nodes with a valuation of $2r$ and two non-terminating nodes with a valuation of at least $2r+1$. When $k\geq r+1$, a non-terminating node at the $k^{th}$ level splits into $p-1$ terminating nodes with a valuation of $k+r$ and one non-terminating node with a valuation of at least $k+1+r$.
    \end{itemize}
\end{theorem}

Before going through the proof of these cases, let us examine some examples which reflect these results.


\begin{example}\label{5n^2+111n-92}
    Consider the polynomial $f(n)=(5n-4)(n+23)=5n^2+111n-92$. This has factors $5n-4$ and $n+23$ with roots $\frac{4}{5}$ and $-23$ respectively. Written as $5$-adic expansions, we have \begin{align*}
        \frac{4}{5}&=4\cdot 5^{-1}+0\cdot 5^0 +0\cdot 5^1+\cdots\\
        -23&=2\cdot 5^0+0\cdot 5^1+4\cdot 5^2 +\cdots,
    \end{align*} which shows us that $\frac{4}{5}\notin \mathbb{Z}_5$ and $-23\in \mathbb{Z}_5$. Hence, the tree for $\{\nu_5(5n^2+111n-92)\}_{n\geq 0}$ will be identical to the tree for $\{\nu_5(n+23)\}_{n\geq 0}$:
    \[
    \xymatrix{
 & & & & & & *+[Fo]{*} \ar@{-}[dlll]_{0} \ar@{-}[dl]_{1} \ar@{-}[d]^{2} \ar@{-}[dr]^{3} \ar@{-}[drrr]^{4} & & & \\
 & & & *+[Fo]{0} & & *+[Fo]{0} & *+[Fo]{*} \ar@{-}[dlll]_{2} \ar@{-}[dl]_{7} \ar@{-}[d]^{12} \ar@{-}[dr]^{17} \ar@{-}[drrr]^{22} & *+[Fo]{0} & & *+[Fo]{0} \\
 & & & *+[Fo]{*} \ar@{-}[dlll]_{2} \ar@{-}[dl]_{27} \ar@{-}[d]^{52} \ar@{-}[dr]^{77} \ar@{-}[drrr]^{102} & & *+[Fo]{1} & *+[Fo]{1} & *+[Fo]{1} & & *+[Fo]{1} \\
 *+[Fo]{2} & & *+[Fo]{2} & *+[Fo]{2} & *+[Fo]{2} & & *+[Fo]{*} \ar@{-}[dlll]_{102} \ar@{-}[dl]_{227} \ar@{-}[d]^{352} \ar@{-}[dr]^{477} \ar@{-}[drrr]^{602} & & & \\
 & & & *+[Fo]{3} & & *+[Fo]{3} & *+[Fo]{3} & *+[Fo]{3} & & *+[Fo]{*} \\
    }
    \] \[\text{Figure 4.2: $5$-Adic valuation tree for $5n^2+111n-92$}\] This tree has one infinite branch corresponding to the $5$-adic expansion for $-23$.
\end{example}

\begin{example}\label{5n^2-82n+32}
      Next, consider $f(n)=(n-16)(5n-2)=5n^2-82n+32$. This has factors $n-16$ and $5n-2$ with roots $16$ and $\frac{2}{5}$ respectively. Written as a $2$-adic expansion, we have \begin{align*}
          16&=0\cdot 2^0+0\cdot 2^1+0\cdot 2^2 +0\cdot 2^3+1\cdot 2^4 +0\cdot 2^5\cdots \\
          \frac{2}{5}&=0\cdot 2^0+1\cdot 2^1+0\cdot 2^2+1\cdot 2^3+1\cdot 2^4+0\cdot 2^5 + \cdots.
      \end{align*} Since $16\equiv \frac{2}{5} \bmod{2^1}$ but $16\not\equiv \frac{2}{5} \bmod{2^2}$, we say $r=1$ and expect the non-terminating node at the $1^{st}$ level to split into two non-terminating nodes. The tree for $\{\nu_2(5n^2-82n+32)\}_{n\geq 0}$, as seen below, has two infinite branches corresponding to the $2$-adic expansions for $16$ and $\frac{2}{5}$.
\[
\xymatrix{
 & & & & & *+[Fo]{*} \ar@{-}[dl]_{0} \ar@{-}[dr]^{1} & & & \\
 & & & & *+[Fo]{*} \ar@{-}[dll]_{0} \ar@{-}[drr]^{2} & & *+[Fo]{0} & & \\
 & & *+[Fo]{*} \ar@{-}[dl]_{0} \ar@{-}[dr]^{4} & & & & *+[Fo]{*} \ar@{-}[dl]_{2} \ar@{-}[dr]^{6} & & \\
 & *+[Fo]{*} \ar@{-}[dl]_{0} \ar@{-}[dr]^{8} & & *+[Fo]{3} & & *+[Fo]{*} \ar@{-}[dl]_{2} \ar@{-}[dr]^{10} & & *+[Fo]{3} & & \\
 *+[Fo]{*} & & *+[Fo]{4} & & *+[Fo]{4} & & *+[Fo]{*} & & & \\
}
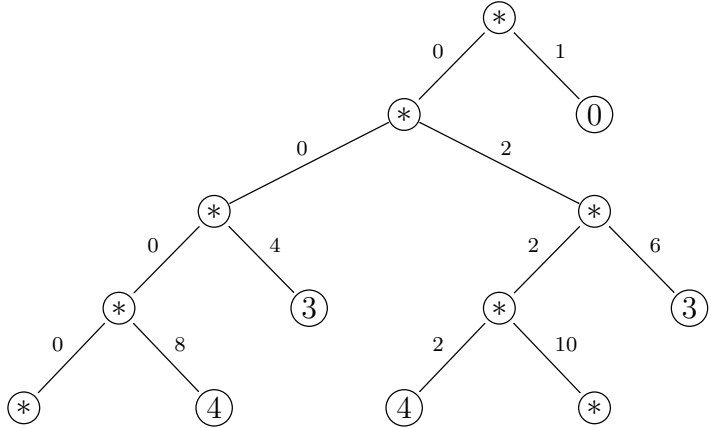
\] \[\text{Figure 4.3: $2$-Adic valuation tree for $5n^2-82n+32$}\] 
\end{example}

Now we will look at the proof for the above theorem. This proof will not be exhaustive due to its similarities to previous proofs and will instead rely upon the properties of $p$-adic valuations.

\vspace{3mm}

\begin{solution}
    (i) Since $-\frac{b_1}{a_1}=-\frac{b_2}{a_2}$ and $\gcd(a_1,b_1)=\gcd(a_2,b_2)=1$, by unique representation, we know that $a_1=a_2$ and $b_1=b_2$. Hence, $a_1n+b_1=a_2n+b_2$ and $f(n)=(a_1n+b_1)^2$. As a result, \[\nu_p(f(n))=\nu_p((a_1n+b_1)^2)=\nu_p(a_1n+b_1)+\nu_p(a_1n+b_1)=2\nu_p(a_1n+b_1).\] If a node on the tree of $a_1n+b_1$ terminates with a valuation of $k$, then $\nu_p(a_1n+p)=k$. Therefore, $\nu_p(f(n))=2k$, which is a fixed valuation. Thus, this same node is terminating on the $p$-adic valuation tree for $f(n)$ and has twice the valuation that $a_1n+b_1$ had. If a node on the tree of $a_1n+b_1$ is non-terminating, then there exist $n_1, n_2$ in this node such that \[\nu_p(a_1n_1+b_1)\neq \nu_p(a_1n_2+b_1).\] Now we see that \begin{align*}
        \nu_p(f(n_1))&=2\nu_p(a_1n_1+b_1)\\
        &\neq 2\nu_p(a_1n_2+b_1)\\
        &=\nu_p(f(n_2))
    \end{align*} and this same node is non-terminating on the tree of $f(n)$.

\vspace{3mm}

    (ii) Since $-\frac{b_2}{a_2}\notin \mathbb{Z}_p$, by Remark \ref{zero point}, $\nu_p(a_2n+b_2)=0$ for all $n\in \mathbb{Z}^{>0}$. Thus, \begin{align*}
        \nu_p(f(n))&=\nu_p((a_1n+b_1)(a_2n+b_2))\\
        &=\nu_p(a_1n+b_1)+\nu_p(a_2n+b_2)\\
        &=\nu_p(a_1n+b_1).
    \end{align*} If a node on the tree of $a_1n+b_1$ terminates with a valuation of $k$, then $\nu_p(a_1n+p)=k$. Therefore, in this same node, $\nu_p(f(n))=k$ is a fixed valuation. Thus, this node is terminating on the $p$-adic valuation tree for $f(n)$ and has the same valuation that $a_1n+b_1$ had. If a node on the tree of $a_1n+b_1$ is non-terminating, then there exist $n_1, n_2$ in this node such that \[\nu_p(a_1n_1+b_1)\neq \nu_p(a_1n_2+b_1).\] Now we see that \begin{align*}
        \nu_p(f(n_1))&=\nu_p(a_1n_1+b_1)\\
        &\neq \nu_p(a_1n_2+b_1)\\
        &=\nu_p(f(n_2))
    \end{align*} and this same node is non-terminating on the tree of $f(n)$.

\vspace{3mm}

    (iii) First, let us consider when $0\leq k\leq r-1$. Since $-\frac{b_1}{a_1}\equiv -\frac{b_2}{a_2} \bmod{p^r}$, we have that $a_1n+b_1\equiv a_2n+b_2 \bmod{p^k}$ for all $0\leq k\leq r$, and the proof is clear from part (i). 

    Next, consider the $r^{th}$ level. Write the two roots as \begin{align*}
        -\frac{b_1}{a_1}&=c_{1,0}+c_{1,1}p+c_{1,2}p^2+\cdots\\
        -\frac{b_2}{a_2}&=c_{2,0}+c_{2,1}p+c_{2,2}p^2+\cdots,
    \end{align*} where $0\leq c_{i,j}\leq p-1$ for $i\in \{1,2\}$ and $j\geq 0$. Thus, $c_{1,j}=c_{2,j}$ when $0\leq j\leq r$ and $c_{1,r+1}\neq c_{2,r+1}$. Using the proof from Proposition \ref{lin.descr}, we can see that on the $r^{th}$ level, the tree for $\{\nu_p(a_1n+b_1)\}_{n\geq 0}$ splits into $p-1$ terminating nodes with a valuation of $r$ and one non-terminating node with a valuation of at least $r+1$. The non-terminating node has the form $c_{1,0}+c_{1,1}p+\cdots+ c_{1,r}p^{r}+np^{r+1}$, where \[-\frac{b_1}{a_1}\equiv c_{1,0}+c_{1,1}p+\cdots+ c_{1,r}p^{r} \bmod{p^{r+1}}.\] The tree for $\{\nu_p(a_2n+b_2)\}_{n\geq 0}$ splits similarly. The non-terminating node has the form $c_{2,0}+c_{2,1}p+\cdots+ c_{2,r}p^{r}+np^{r+1}$, where \[-\frac{b_2}{a_2}\equiv c_{2,0}+c_{2,1}p+\cdots+ c_{2,r}p^{r} \bmod{p^{r+1}}.\] Thus, we have \begin{align*}
        -\frac{b_1}{a_1}&\not\equiv -\frac{b_2}{a_2} \bmod{p^{r+1}}\\
        c_{1,0}+c_{1,1}p+\cdots+ c_{1,r}p^{r}&\not\equiv c_{2,0}+c_{2,1}p+\cdots+ c_{2,r}p^{r} \bmod{p^{r+1}},
    \end{align*} and these two nodes are not the same. 

    Without loss of generality, we will show that one of these nodes is non-terminating. To see that this branch has a non-terminating node, consider \[c_{1,0}+c_{1,1}p+c_{1,2}p^2+\cdots+c_{1,r+1}p^{r+1}+\cdots +c_{1,s}p^s\] with $s$ large enough. This natural number is within the node and is arbitrarily close to the root $-\frac{b_1}{a_1}$. Since $f(n)$ is a continuous function, $f(c_{1,0}+c_{1,1}p+c_{1,2}p^2+\cdots+c_{1,r+1}p^{r+1}+\cdots +c_{1,s}p^s)$ is arbitrarily close to $f\left(-\frac{b_1}{a_1}\right)=0$. Thus, \[\nu_p(f(c_{1,0}+c_{1,1}p+c_{1,2}p^2+\cdots+c_{1,r+1}p^{r+1}+\cdots +c_{1,s}p^s))\] is arbitrarily large. Now, consider the natural number \[c_{1,0}+c_{1,1}p+c_{1,2}p^2+\cdots+(c_{1,r+1}+\alpha)p^{r+1}+\cdots +c_{1,s}p^s\not\equiv -\frac{b_2}{a_2} \bmod{p^{s+1}}\] for some $\alpha\in \mathbb{Z}$. This is also within the same node, but not as close to $-\frac{b_1}{a_1}$. This number is also not close to $-\frac{b_2}{a_2}$. Again, since $f(n)$ is a continuous function with two roots, $f(c_{1,0}+c_{1,1}p+c_{1,2}p^2+\cdots+(c_{1,r+1}+\alpha)p^{r+1}+\cdots +c_{1,s}p^s)$ is not close to $f\left(-\frac{b_1}{a_1}\right)=f\left(-\frac{b_2}{a_2}\right)=0$ so \[\nu_p(f(c_{1,0}+c_{1,1}p+c_{1,2}p^2+\cdots+(c_{1,r+1}+\alpha)p^{r+1}+\cdots +c_{1,s}p^s))\] is some fixed finite value which is different than the previous valuation. Since this node has multiple different valuations, it is non-terminating.
    
    This results in two distinct non-terminating nodes, one derived from each of the roots of the polynomial. The other $p-2$ nodes are terminating for both the tree of $\{\nu_p(a_1n+b_1)\}_{n\geq 0}$ and $\{\nu_p(a_2n+b_2)\}_{n\geq 0}$, so their valuations are fixed, specifically with a valuation of $r$. Since $\nu_p(f(n))=\nu_p(a_1n+b_1)+\nu_p(a_2n+b_2)$, the valuations at these nodes for the tree of $\{\nu_p(f(n))\}_{n\geq 0}$ are fixed. Thus, these nodes are terminating and have a valuation of $2r$.

    Lastly, consider the $k^{th}$ level when $k\geq r+1$. Without loss of generality, consider a non-terminating node that splits into $p$ branches of the form $c_{1,0}+c_{1,1}p+\cdots+c_{1,k-1}p^{k-1}+i_kp^k+np^{k+1}$, where $0\leq i_k\leq p-1$ and \[-\frac{b_1}{a_1}\equiv c_{1,0}+c_{1,1}p+\cdots+c_{1,k-1}p^{k-1}+i_kp^k \bmod{p^k}.\] Since $-\frac{b_1}{a_1}\not\equiv -\frac{b_2}{a_2} \bmod{p^{r+1}}$, Proposition \ref{lin.descr} tells us that within the tree of $\{\nu_p(a_2n+b_2)\}_{n\geq 0}$, the nodes the form $c_{1,0}+c_{1,1}p+\cdots+c_{1,k-1}p^{k-1}+i_kp^k+np^{k+1}$ are contained within the node of the form $c_{1,0}+c_{1,1}p+\cdots+c_{1,r+1}p^{r+1}$, which terminates with a valuation of $r$. Within the tree of $\{\nu_p(a_1n+b_1)\}_{n\geq 0}$, the non-terminating node on the $k^{th}$ level splits into $p-1$ terminating nodes with a valuation of $k$ and one non-terminating node with a valuation of at least $k+1$. For the nodes that terminate on the tree of $\{\nu_p(a_1n+b_1)\}_{n\geq 0}$, we have \begin{align*}
        \nu_p(f(n))&=\nu_p(a_1n+b_1)+\nu_p(a_2n+b_2)\\
        &=k+r,
    \end{align*} which is a fixed valuation. Hence, $p-1$ nodes terminate with a valuation of $k+r$. For the node that is non-terminating on the tree of $\{\nu_p(a_1n+b_1)\}_{n\geq 0}$, we have \begin{align*}
        \nu_p(f(n))&=\nu_p(a_1n+b_1)+\nu_p(a_2n+b_2)\\
        &\geq k+1+r,
    \end{align*} which is not a fixed valuation since $\nu_p(a_1n+b_1)$ is not fixed, while $\nu_p(a_2n+b_2)$ is fixed. Hence, this node is non-terminating with a valuation of at least $k+1+r$.
\end{solution}

\vspace{3mm}

We see that the $p$-adic valuation tree for the product of two linear factors has an infinite branch for each distinct $p$-adic integer root. Again, the infinite branch corresponds to the $p$-adic expansion of that root.

Using Lemma \ref{coef}, we can extend this theorem to describe the tree for any factorable quadratic polynomial with coefficients in $\mathbb{Q}$.

\begin{example}\label{3/2(n-9)(3n+1)}
    Consider the polynomial $f(n)=\frac{3}{2}(n-9)(3n+1)=\frac{9}{2}n^2-39n-\frac{27}{2}$. This has factors $n-9$ and $3n+1$ with roots $9$ and $-\frac{1}{3}$ respectively. Written as $2$-adic expansions, we have \begin{align*}
        9&=1\cdot 2^0+0\cdot 2^1+0\cdot 2^2+1\cdot 2^3 \cdots \\
        -\frac{1}{3}&=1\cdot 2^0+0\cdot 2^1+1\cdot 2^2+0\cdot 2^3\cdots.
    \end{align*} Since $9\equiv -\frac{1}{3} \bmod{2^2}$, but $9\not\equiv -\frac{1}{3} \bmod{2^3}$, we say $r=2$ and expect our two infinite branches to separate at the $2^{nd}$ level. Furthermore, the factor of $\frac{3}{2}$ means that the valuations of the terminating nodes change by $\nu_2\left(\frac{3}{2}\right)=-1$. The tree for $\{\nu_2\left(\frac{9}{2}n^2-39n-\frac{27}{2}\right)\}_{n\geq 0}$ matches our expected behavior with the two infinite branches and valuations one less due to the factor of $\frac{3}{2}$.
\[
\xymatrix{
 & & & *+[Fo]{*} \ar@{-}[dl]_{0} \ar@{-}[dr]^{1} & & & \\
 & & *+[Fo]{-1} & & *+[Fo]{*} \ar@{-}[dl]_{1} \ar@{-}[dr]^{3} & & \\
 & & & *+[Fo]{*} \ar@{-}[dll]_{1} \ar@{-}[drr]^{5} & & *+[Fo]{1} & \\
 & *+[Fo]{*} \ar@{-}[dl]_{1} \ar@{-}[dr]^{9} & & & & *+[Fo]{*} \ar@{-}[dl]_{5} \ar@{-}[dr]^{13} & \\
 *+[Fo]{4} & & *+[Fo]{*} & & *+[Fo]{*} & & *+[Fo]{4} \\
}
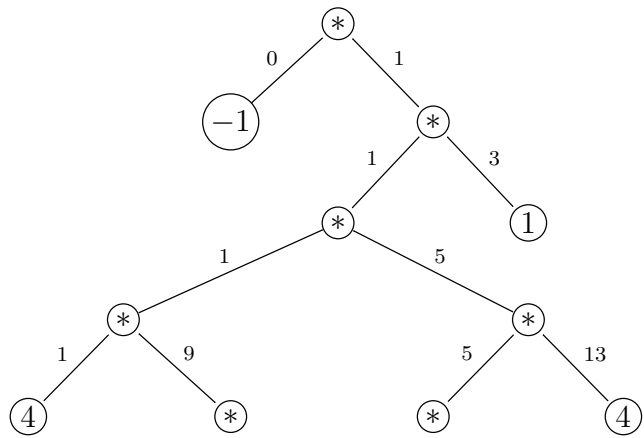
\] \[\text{Figure 4.4: $2$-Adic valuation tree for $\frac{9}{2}n^2-39n-\frac{27}{2}$}\] 
\end{example}

\section{Combining three or more linear trees}

We can extend our previous results to describe the $p$-adic valuation tree for a polynomial comprised of any finite number of linear factors.

\begin{lemma}
    Let $f(n)=(a_1n+b_1)(a_2n+b_2)\cdots (a_rn+b_r)$, where $a_i,b_i\in \mathbb{Z}$, $a_i\neq 0$, and $\gcd(a_i,b_i)=1$ for all $1\leq i\leq r$. Suppose the root of one of the factors, $-\frac{b_j}{a_j}$, is not a $p$-adic integer. Then the tree of $\{\nu_p(f(n))\}_{n\geq 0}$ is identical to the tree of $\Big\{\nu_p\left(\frac{f(n)}{a_jn+b_j}\right)\Big\}_{n\geq 0}$.
\end{lemma}

\begin{solution}
    Suppose $-\frac{b_j}{a_j}\notin \mathbb{Z}_p$ for some $1\leq j\leq r$. Then, by Remark \ref{zero point}, we know $\nu_p(a_jn+b_j)=0$. Therefore, we have \begin{align*}
        \nu_p(f(n))&= \nu_p(f(n))-\nu_p(a_jn+b_j)+\nu_p(a_jn+b_j)\\
        &=\nu_p\left(\frac{f(n)}{a_jn+b_j}\right)+\nu_p(a_jn+b_j)\\
        &=\nu_p\left(\frac{f(n)}{a_jn+b_j}\right).
    \end{align*} 
\end{solution}

\begin{theorem}\label{many lin}
   Let $f(n)=(a_1n+b_1)(a_2n+b_2)\cdots (a_rn+b_r)$, where $a_i,b_i\in \mathbb{Z}$, $a_i\neq 0$, and $\gcd(a_i,b_i)=1$ for all $1\leq i\leq r$. Suppose the root of each factor, $-\frac{b_i}{a_i}$, is a $p$-adic integer. Write their $p$-adic expansions as \begin{align*}
       -\frac{b_1}{a_1}&=c_{1,0}+c_{1,1}p+c_{1,2}p^2+\cdots, \\
       -\frac{b_2}{a_2}&=c_{2,0}+c_{2,1}p+c_{2,2}p^2+\cdots, \\
       &\hspace{2mm} \vdots \\
       -\frac{b_r}{a_r}&=c_{r,0}+c_{r,1}p+c_{r,2}p^2+\cdots,
   \end{align*} where $0\leq c_{i,j}\leq p-1$ for all $1\leq i\leq r$ and $j\geq 0$. For any $k\geq 0$, consider a node on the $k^{th}$ level. Write the terms in this node as $d_0+d_1p+d_2p^2+\cdots +d_{k-1}p^{k-1}+np^k$, where $0\leq d_i\leq p-1$. For all $1\leq i\leq r$, find the smallest $j< k$ such that $d_j\neq c_{i,j}$ and denote this $j$ as $v_i$. If, for any $1\leq i\leq r$, no such $j$ exists, then the node is non-terminating. Otherwise, the node is terminating with valuation \[\sum_{i=1}^r v_i.\]
\end{theorem}

\begin{solution}
    First, suppose there exists some root \[-\frac{b_i}{a_i}=c_{i,0}+c_{i,1}p+c_{i,2}p^2+\cdots+c_{i,k}p^k +\cdots,\] where $d_j=c_{i,j}$ for all $j<k$. Then the node of the form $d_0+d_1p+d_2p^2+\cdots +d_{k-1}p^{k-1}+np^k$ contains the truncation of this root \[c_{i,0}+c_{i,1}p+c_{i,2}p^2+\cdots+c_{i,k}p^k +\cdots+c_{i,s}p^s,\] where $s$ is arbitrarily large. We can show this is contained in the node by setting \[n=c_{i,k}+c_{i,k+1}p+\cdots +c_{i,s}p^{s-k}.\] To see that this same node is non-terminating in the $p$-adic valuation tree of $f(n)$, notice that the natural number \[c_{i,0}+c_{i,1}p+c_{i,2}p^2+\cdots+c_{i,k}p^k +\cdots+c_{i,s}p^s\] for large enough $s$ is within the node and is arbitrarily close to the root $-\frac{b_i}{a_i}$. Since $f(n)$ is a continuous function, $f(c_{i,0}+c_{i,1}p+c_{i,2}p^2+\cdots+c_{i,k}p^k +\cdots+c_{i,s}p^s)$ is arbitrarily close to $f\left(-\frac{b_i}{a_i}\right)=0$. Thus, \[\nu_p(f(c_{i,0}+c_{i,1}p+c_{i,2}p^2+\cdots+c_{i,k}p^k +\cdots+c_{i,s}p^s))\] is arbitrarily large. Now consider the natural number \[c_{i,0}+c_{i,1}p+c_{i,2}p^2+\cdots+(c_{i,k}+\alpha)p^k +\cdots+c_{i,s}p^s\] where $\alpha\in \mathbb{Z}$ is chosen such that this natural number does not match a truncation of $-\frac{b_j}{a_j}$ for any $1\leq j\leq r$. This is also within the same node, but not as close to $-\frac{b_i}{a_i}$ or any other $-\frac{b_j}{a_j}$. Again, since $f(n)$ is continuous and has $r$ roots, $f(c_{i,0}+c_{i,1}p+c_{i,2}p^2+\cdots+(c_{i,k}+\alpha)p^k +\cdots+c_{i,s}p^s)$ is not as close to $f\left(-\frac{b_j}{a_j}\right)=0$, so \[\nu_p(f(c_{i,0}+c_{i,1}p+c_{i,2}p^2+\cdots+(c_{i,k}+\alpha)p^k +\cdots+c_{i,s}p^s))\] is some fixed finite value that is different than the previous valuation. Hence, this node has two terms with different valuations and is therefore non-terminating.

    Next, suppose $v_i$ exists for all $1\leq i\leq r$. Hence, there exists $j<k$ such that $d_j\neq c_{i,j}$. Choose some $1\leq t\leq r$ and set $0\leq j_t<k$ to be the smallest integer, where $d_{j_t}\neq c_{t,j_t}$. Here, $v_t=j_t$. Then, the node of the form \[d_0+d_1p+d_2p^2+\cdots +d_{j_t-1}p^{j_t-1}+np^{j_t}\] contains the truncation of the root $-\frac{b_t}{a_t}$, and is therefore non-terminating. By Proposition \ref{lin.descr}, we know that on the tree for $\{\nu_p(a_tn+b_t)\}_{n\geq 0}$, this node splits into one non-terminating node and $p-1$ terminating nodes with a valuation of $j_t$. The node of the form $d_0+d_1p+d_2p^2+\cdots +d_{j_t}p^{j_t}+np^{j_t+1}$ contains the node of the form $d_0+d_1p+d_2p^2+\cdots +d_{k-1}p^{k-1}+np^k$ but does not contain the truncation of the root $-\frac{b_t}{a_t}$. Hence, this node is terminating with a valuation of $j_t$ on the tree for $\{\nu_p(a_tn+b_t)\}_{n\geq 0}$. Now, looking at the node of the form \[d_0+d_1p+d_2p^2+\cdots +d_{k-1}p^{k-1}+np^k,\] we know it must also have a valuation of $j_t$ for the tree for $\{\nu_p(a_tn+b_t)\}_{n\geq 0}$ since it is contained within a terminating node. 

    We can repeat this process for all $1\leq i\leq r$ and retrieve $r$ fixed valuations, one for each factor. Hence, at the node of the form \[d_0+d_1p+d_2p^2+\cdots +d_{k-1}p^{k-1}+np^k,\] we have \begin{align*}
        \nu_p(f(n))&= \nu_p(a_1n+b_1)+\nu_p(a_2n+b_2)+\cdots + \nu_p(a_rn+b_r)\\
        &=v_1+v_2+\cdots +v_r\\
        &=\sum_{i=1}^r v_i,
    \end{align*} which is a fixed valuation, signifying a terminating node.
\end{solution}

\vspace{3mm}

This is a different method of constructing trees compared to the previous results. Instead of describing how the non-terminating nodes split within the tree, this proposition tells us the status of individual nodes. Regardless, we can still discuss some general properties of the overall tree. 

The following corollary relates to the result from Medina et al. \cite{medina} which showed that a $p$-adic valuation sequence generated by a polynomial with integer coefficients is unbounded if and only if it has a $p$-adic integer root. The unbounded behavior appears as an infinite branch and, in fact, this corollary will show the exact number of distinct $p$-adic integer roots.

\begin{corollary}
     Let $f(n)=(a_1n+b_1)(a_2n+b_2)\cdots (a_rn+b_r)$, where $a_i,b_i\in \mathbb{Z}$, $a_i\neq 0$, and $\gcd(a_i,b_i)=1$ for all $1\leq i\leq r$. Then, there are a number of infinite branches equal to the number of distinct $p$-adic integer roots. Moreover, the form of each infinite branch corresponds to the $p$-adic expansion for one of these roots.
\end{corollary}

\begin{solution}
  From Proposition \ref{lin.descr}, we know that the tree $\{\nu_p(a_in+b_i)\}_{n\geq 0}$ has an infinite branch if and only if $-\frac{b_i}{a_i}$ is a $p$-adic integer. We also know that this infinite branch corresponds to the $p$-adic expansion for this root. Using a similar process as seen in the proof for Theorem \ref{many lin}, we can show that every non-terminating node on the infinite branch for the tree $\{\nu_p(a_in+b_i)\}_{n\geq 0}$ will also be non-terminating on the tree $\{\nu_p(f(n))\}_{n\geq 0}$. Thus, this infinite branch remains on the combined tree and still corresponds to the $p$-adic expansion for that root. This is done for every factor, creating an infinite branch for each $p$-adic integer root whose structure matches the $p$-adic expansion of this root. Repeated roots result in the same infinite branch, so only distinct $p$-adic integer roots create distinct infinite branches.
\end{solution}

\vspace{3mm}

Thus, when constructing these trees, we can easily create the infinite branches by looking at the $p$-adic expansions of the roots. In the case of linear factors, all other nodes will be terminating. 

\begin{example}\label{n^3-8n^2-60n+144}
    Consider the polynomial \[f(n)=(n-2)(n+6)(n-12)=n^3-8n^2-60n+144.\] This has factors $n-2$, $n+6$, and $n-12$ with roots $2$, $-6$, and $12$ respectively. Written as $2$-adic expansions, we have \begin{align*}
        2&=0\cdot 2^0+1\cdot 2^1+0\cdot 2^2+0\cdot 2^3 \cdots \\
        -6&=0\cdot 2^0+1\cdot 2^1+0\cdot 2^2+1\cdot 2^3 \cdots \\
        12&=0\cdot 2^0+0\cdot 2^1+1\cdot 2^2+1\cdot 2^3 \cdots.
    \end{align*} We expect there to be three infinite branches, each resembling the infinite branches on their individual trees. We see the tree for $\{\nu_2(n^3-8n^2-60n+144)\}_{n\geq 0}$ below.
\[
\xymatrix{
 & & & & & *+[Fo]{*} \ar@{-}[dl]_{0} \ar@{-}[dr]^{1} & & & \\
 & & & & *+[Fo]{*} \ar@{-}[dlll]_{0} \ar@{-}[drrr]^{2} & & *+[Fo]{0} & & \\
 & *+[Fo]{*} \ar@{-}[dl]_{0} \ar@{-}[dr]^{4} & & & & & & *+[Fo]{*} \ar@{-}[dl]_{2} \ar@{-}[dr]^{6} & & \\
 *+[Fo]{4} & & *+[Fo]{*} \ar@{-}[dl]_{4} \ar@{-}[dr]^{12} & & & & *+[Fo]{*} \ar@{-}[dl]_{2} \ar@{-}[dr]^{10} & & *+[Fo]{5} & \\
 & *+[Fo]{5} & & *+[Fo]{*} & & *+[Fo]{*} & & *+[Fo]{*} \\
}
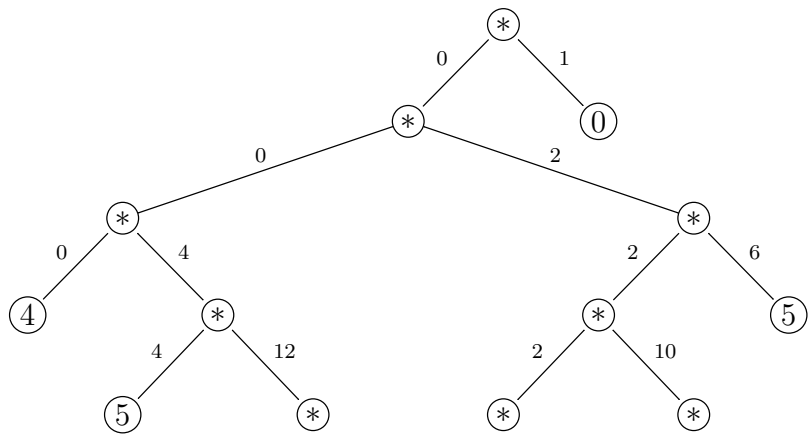
\] \[\text{Figure 5.1: $2$-Adic valuation tree for $n^3-8n^2-60n+144$}\]
\end{example}

Notice that each of the infinite branches on the tree corresponds to one of the distinct $p$-adic integer roots of the polynomial.

Again, using Lemma \ref{coef}, we can extend this theorem to describe the tree for any completely factorable polynomial of any degree with coefficients in $\mathbb{Q}$.

\begin{example}\label{3/5(n+6)(2n-5)(3n-4)}
    Consider $p=3$ and the polynomial \[f(n)=\frac{3}{5}(n+6)(2n-5)(3n-4)=\frac{18}{5}n^3+\frac{39}{5}n^2-\frac{352}{5}n+72.\] This has factors $n+6$, $2n-5$, and $3n-4$ with roots $-6$, $\frac{5}{2}$, and $\frac{4}{3}$ respectively. Written as $3$-adic expansions, we have \begin{align*}
         -6&=0\cdot 3^0+1\cdot 3^1+2\cdot 3^2+2\cdot 3^3+2\cdot 3^4+\cdots \\
            \frac{5}{2}&= 1\cdot 3^0+2\cdot 3^1+1\cdot 3^2+1\cdot 3^3+1\cdot 3^4+\cdots \\
            \frac{4}{3}&=1 \cdot 3^{-1} +1\cdot 3^0 +0\cdot 3^1 +0\cdot 3^2+0\cdot 3^3+0\cdot 3^4+\cdots.
    \end{align*} Since $\frac{4}{3}$ is not a $3$-adic integer, we expect there to be two infinite branches corresponding to the $3$-adic expansions for the roots $-6$ and $\frac{5}{2}$. Also, the valuations of the terminating node will be one greater since $\nu_3\left(\frac{3}{5}\right)=1$. The tree for $\{\nu_3\left(\frac{18}{5}n^3+\frac{39}{5}n^2-\frac{352}{5}n+72\right)\}_{n\geq 0}$ matches our expectations.
  \[
\xymatrix{
 & & & & & *+[Fo]{*} \ar@{-}[dllll]_{0} \ar@{-}[d]^{1} \ar@{-}[dr]^{2} & & \\
 & *+[Fo]{*}  \ar@{-}[dl]_{0} \ar@{-}[d]^{3} \ar@{-}[dr]^{6} & & & & *+[Fo]{*} \ar@{-}[dl]_{1} \ar@{-}[d]^{4} \ar@{-}[dr]^{7} & *+[Fo]{1} & \\
 *+[Fo]{2} & *+[Fo]{*} \ar@{-}[dl]_{3} \ar@{-}[d]^{12} \ar@{-}[dr]^{21} & *+[Fo]{2} & & *+[Fo]{2} & *+[Fo]{2} & *+[Fo]{*} \ar@{-}[dl]_{7} \ar@{-}[d]^{16} \ar@{-}[dr]^{25} & \\
 *+[Fo]{3} & *+[Fo]{3} & *+[Fo]{*} \ar@{-}[dl]_{21} \ar@{-}[d]^{48} \ar@{-}[dr]^{75} & & & *+[Fo]{3} & *+[Fo]{*} \ar@{-}[dl]_{16} \ar@{-}[d]^{43} \ar@{-}[dr]^{70} & *+[Fo]{3} \\
 & *+[Fo]{4} & *+[Fo]{4} & *+[Fo]{*}  \ar@{-}[dl]_{75} \ar@{-}[d]^{156} \ar@{-}[dr]^{237} & & *+[Fo]{4} & *+[Fo]{*} \ar@{-}[dl]_{43} \ar@{-}[d]^{124} \ar@{-}[dr]^{205} & *+[Fo]{4} \\
 & & *+[Fo]{5} & *+[Fo]{5} & *+[Fo]{*} & *+[Fo]{5} & *+[Fo]{*} & *+[Fo]{5}
}
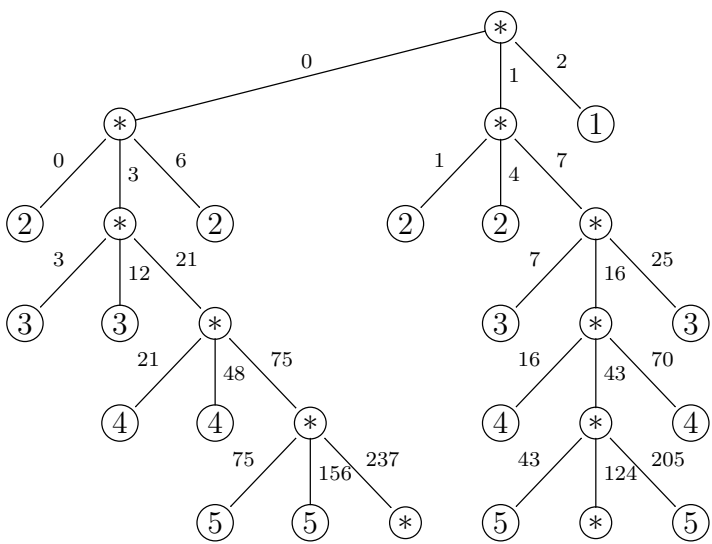
\] \[\text{Figure 5.2: $3$-Adic valuation tree for $\frac{18}{5}n^3+\frac{39}{5}n^2-\frac{352}{5}n+72$}\] 
\end{example}

\section{Combining higher-degree trees}

Now, let us discuss functions that are products of higher-degree polynomials. The descriptions for their trees are far less complete and an area for further research.

First, we can state that infinite branches from the factors remain in the combined tree. 

\begin{proposition}
    Let $f(n)=g_1(n)g_2(n)\cdots g_r(n)$, where $g_i(n)\in \mathbb{Z}[n]$ for all $1\leq i\leq r$. Suppose $g_j(n)$ has a root in $\mathbb{Z}_p$ for some prime $p$. Write this root as its $p$-adic expansion \[c_0+c_1p+c_2p^2+\cdots\] where $0\leq c_k\leq p-1$ for all $k\geq 0$. Then, on the tree for $\{\nu_p(f(n))\}_{n\geq 0}$, the nodes of the form $c_0+c_1+\cdots +c_tp^t+np^{t+1}$ are non-terminating for all $t\geq 0$.
\end{proposition}

From Theorem \ref{almo}, we know that when a function has a $p$-adic integer root, it results in an infinite branch on its $p$-adic valuation tree. Moreover, this infinite branch displays the coefficients of the root's $p$-adic expansion. For a product, the root of one factor is also a root for the product, resulting in the same infinite branch on our new tree.

\vspace{3mm}

\begin{solution}
    Consider the factor $g_j(n)$ and a root of this polynomial, $\alpha$, with $p$-adic expansion \[\alpha=c_0+c_1p+c_2p^2+\cdots,\] where $0\leq c_k\leq p-1$ for all $k\geq 0$. Then, $f(\alpha)=0$, and we see that $\alpha$ is also a root for the product. 
    
    Now consider a node on the $t^{th}$ level of the form $c_0+c_1+\cdots +c_tp^t+np^{t+1}$ for any $t\geq 0$. To see that this node is non-terminating on the $p$-adic valuation tree of $f(n)$, consider \[c_0+c_1p+c_2p^2+\cdots+c_tp^t+\cdots +c_sp^s\] with $s$ large enough. This natural number is within the node and is arbitrarily close to the root $\alpha$. Since $f(n)$ is a continuous function, $f(c_0+c_1p+c_2p^2+\cdots+c_tp^t+\cdots +c_sp^s)$ is arbitrarily close to $f(\alpha)=0$. Thus, \[\nu_p(f(c_0+c_1p+c_2p^2+\cdots+c_tp^t+\cdots +c_sp^s))\] is arbitrarily large. Now consider the natural number \[c_0+c_1p+c_2p^2+\cdots+c_tp^t+(c_{t+1}+\beta)p^{t+1}+\cdots +c_sp^s,\] where $\beta\in \mathbb{Z}$ is chosen such that this natural number is not another root of $g_j(n)$ nor a root for any other $g_i(n)$. This is also within the same node but not as close to $\alpha$ or any root of $f(n)$. Again, since $f(n)$ is continuous and has a finite number of roots, $f(c_0+c_1p+c_2p^2+\cdots+c_tp^t+(c_{t+1}+\beta)p^{t+1}+\cdots +c_sp^s)$ is not as close to $f(\alpha)=0$, so \[\nu_p(f(c_0+c_1p+c_2p^2+\cdots+c_tp^t+(c_{t+1}+\beta)p^{t+1}+\cdots +c_sp^s))\] is some fixed finite value that is different than the previous valuation. Hence, this node has two terms with different valuations and is therefore non-terminating.
\end{solution}

\vspace{3mm}

We can also show some nodes which are definitively terminating on the $p$-adic valuation tree.

\begin{proposition}
    Let $f(n)=g_1(n)g_2(n)\cdots g_r(n)$, where $g_i(n)\in \mathbb{Z}[n]$ for all $1\leq i\leq r$ as above. Consider a node on the $p$-adic valuation tree. If that node, or a node at a higher level that contains the aforementioned node, is terminating on the tree for every $\{\nu_p(g_i(n))\}_{n\geq 0}$, then the node is terminating on the tree for $\{\nu_p(f(n))\}_{n\geq 0}.$ Moreover, the valuation of this node is the summation of the valuations for the individual factors at that node.
\end{proposition}

\begin{solution}
    Consider a node on the $t^{th}$ level of the form $c_0+c_1+\cdots +c_tp^t+np^{t+1}$ for any $t\geq 0$. Suppose that this node is terminating, or a node that contains it terminates for every tree $\{\nu_p(g_i(n))\}_{n\geq 0}$. Then every tree has a fixed valuation at that node. Let us define $v_i$ to be the valuation at this node for the tree $\{\nu_p(g_i(n))\}_{n\geq 0}$. Then, at this node, we have \begin{align*}
        \nu_p(f(c_0+c_1+\cdots +c_tp^t+np^{t+1}))&= \nu_p(g_1(c_0+c_1+\cdots +c_tp^t+np^{t+1}))\\
        &+\nu_p(g_2(c_0+c_1+\cdots +c_tp^t+np^{t+1}))\\
        &+\cdots +\nu_p(g_r(c_0+c_1+\cdots +c_tp^t+np^{t+1}))\\
        &= v_1+v_2+\cdots +v_r,
    \end{align*} which is a fixed value. Hence, this node is terminating for the tree $\{\nu_p(f(n))\}_{n\geq 0}$, and its valuation is the sum of the valuation of the factors.
\end{solution}

\vspace{3mm}

Let us look at an example of this. While we will only use $p=2$ for simplicity, be aware that this works for any prime $p$.

\begin{example}\label{n^4+11n^2+28}
    Consider the polynomial $f(n)=(n^2+4)(n^2+7)=n^4+11n^2+28$. The factor $n^2+7$ has roots $\pm \sqrt{-7}$ with their $2$-adic expansions being
    \begin{align*}
        \sqrt{-7}&=1\cdot 2^0+ 0\cdot 2^1+ 1\cdot 2^2 +0\cdot 2^3+\cdots\\
        -\sqrt{-7}&=1\cdot 2^0+ 1\cdot 2^1+ 0\cdot 2^2+ 1\cdot 2^3+\cdots.
    \end{align*} The factor $n^2+4$ has roots $\pm \sqrt{-4}$, which are not $2$-adic integers. 
    
    The tree for $\{\nu_2(n^2+7)\}_{n\geq 0}$ has two infinite branches. 
    \[
    \xymatrix{
 & & *+[Fo]{*} \ar@{-}[dl]_{0} \ar@{-}[dr]^{1} & & & & & \\
 & *+[Fo]{0} & & *+[Fo]{*} \ar@{-}[dll]_{1} \ar@{-}[drrr]^{3} & & & & \\ 
 & *+[Fo]{*} \ar@{-}[dl]_{1} \ar@{-}[dr]^{5} & & & & & *+[Fo]{*} \ar@{-}[dl]_{3} \ar@{-}[dr]^{7} & \\
 *+[Fo]{3} & & *+[Fo]{*} \ar@{-}[dl]_{5} \ar@{-}[dr]^{13} & & & *+[Fo]{*} \ar@{-}[dl]_{3} \ar@{-}[dr]^{11} & & *+[Fo]{3} \\
 & *+[Fo]{*} & & *+[Fo]{4} & *+[Fo]{4} & & *+[Fo]{*} & \\
    }
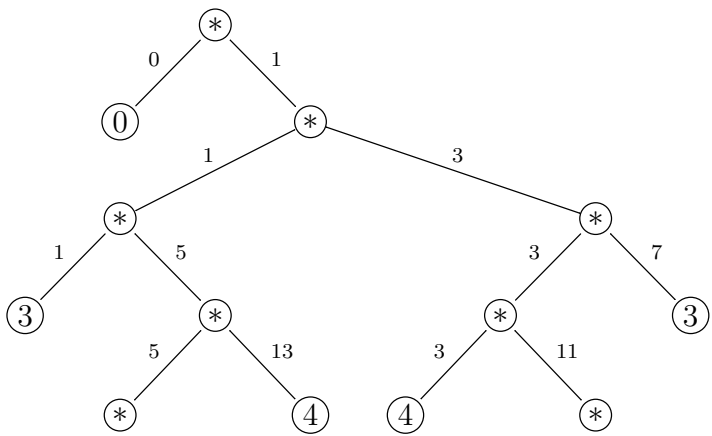
\] \[\text{Figure 6.1: $2$-Adic valuation tree of $n^2+7$}\]
    The tree for $\{\nu_2(n^2+4)\}_{n\geq 0}$ has no infinite branches.
    
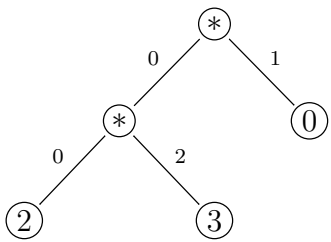
\[
    \xymatrix{
    & & *+[Fo]{*} \ar@{-}[dl]_{0} \ar@{-}[dr]^{1} & \\
    & *+[Fo]{*} \ar@{-}[dl]_{0} \ar@{-}[dr]^{2} & & *+[Fo]{0} \\
    *+[Fo]{2} & & *+[Fo]{3} & \\
    }
\] \[\text{Figure 6.2: $2$-Adic valuation tree of $n^2+4$}\]

    Thus, we expect the tree for $\{\nu_2(n^4+11n^2+28)\}_{n\geq 0}$ to have two infinite branches which match the infinite branches from the tree $\{\nu_2(n^2+7)\}_{n\geq 0}$. We also expect that the nodes of the form $2^2n$, $2^2n+2$, $2^3n+1$, $2^3n+2^2+2^1+1$, $2^4n+2^2+1$, and $2^4n+2^1+1$ to be terminating with valuations $2$, $3$, $3$, $3$, $4$, and $4$ respectively.

    The tree for $\{\nu_2(n^4+11n^2+28)\}_{n\geq 0}$ confirms our hypothesis.
    \[
    \xymatrix{
    & & & *+[Fo]{*} \ar@{-}[dll]_{0} \ar@{-}[drr]^{1} & & & & & & \\
    & *+[Fo]{*} \ar@{-}[dl]_{0} \ar@{-}[dr]^{2} & & & & *+[Fo]{*} \ar@{-}[dll]_{1} \ar@{-}[drrr]^{3} & & & & \\   
    *+[Fo]{2} & & *+[Fo]{3} & *+[Fo]{*} \ar@{-}[dl]_{1} \ar@{-}[dr]^{5} & & & & & *+[Fo]{*} \ar@{-}[dl]_{3} \ar@{-}[dr]^{7} & \\
    & & *+[Fo]{3} & & *+[Fo]{*} \ar@{-}[dl]_{5} \ar@{-}[dr]^{13} & & & *+[Fo]{*} \ar@{-}[dl]_{3} \ar@{-}[dr]^{11} & & *+[Fo]{3} \\
    & & & *+[Fo]{*} & & *+[Fo]{4} & *+[Fo]{4} & & *+[Fo]{*} & \\
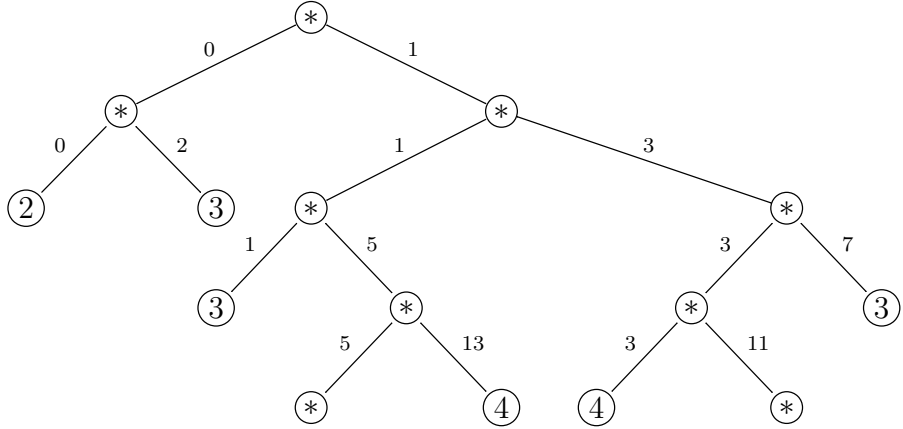
    }\] \[\text{Figure 6.3: $2$-Adic valuation tree of $n^4+11n^2+28$}\]
\end{example}

However, as we may already see, this does not cover every node on these trees.

\begin{remark}
     The trees for higher degree polynomials may contain nodes that are non-terminating, but are not part of an infinite branch. In that case, that same node on the combined tree for a product of higher degree polynomials may not necessarily be non-terminating.
\end{remark}

\begin{example}\label{bad}
    First, notice that the tree for $\{\nu_2(n^2+1)\}_{n\geq 0}$ has a non-terminating node that is not part of an infinite branch. This is clear, considering there is no infinite branch on this tree. 
    \[
\xymatrix{
 & *+[Fo]{*} \ar@{-}[dl]_{0} \ar@{-}[dr]^{1} & \\
 *+[Fo]{0} & & *+[Fo]{1} 
}
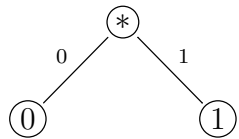
\] \[\text{Figure 6.4: $2$-Adic valuation tree of $n^2+1$}\] 

Also, notice that the tree for $\{\nu_2(n^2+2)\}_{n\geq 0}$ has a non-terminating node that is not part of an infinite branch.  
    \[
\xymatrix{
 & *+[Fo]{*} \ar@{-}[dl]_{0} \ar@{-}[dr]^{1} & \\
 *+[Fo]{1} & & *+[Fo]{0} 
}
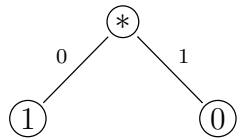
\] \[\text{Figure 6.5: $2$-Adic valuation tree of $n^2+2$}\] 

Their product is $f(n)=(n^2+1)(n^2+2)=n^4+3n^2+2$. The combined tree for $\{\nu_2(n^4+3n^2+2)\}_{n\geq 0}$ is terminating at the node that was non-terminating in the individual trees.
\[
\xymatrix{
 *+[Fo]{1}
}
\] \[\text{Figure 6.6: $2$-Adic valuation tree of $n^4+3n^2+2$}\]
\end{example}

There may be additional ways to describe $p$-adic valuation trees based on their factors. As seen above in Example \ref{bad}, the non-terminating node in the individual trees disappears in the combined tree. This is due to the split branches of that node having reversed valuations in the two factors, specifically, either $0$ or $1$. In terms of the sequence, this means that $n^2+1$ alternates from $0$ to $1$ while $n^2+2$ alternates from $1$ to $0$. Thus, when we multiply them together, the corresponding terms in the sequences add together, so each term has a $0$ and $1$, resulting in a constant valuation. There is further investigation that can be done to discuss exactly when this phenomenon occurs.

\section{Conclusion}

Using the results detailed above, we can describe many properties of $p$-adic valuation trees generated by polynomial sequences. If the polynomial is factorable into linear parts, we can describe the tree in its entirety, identifying the terminating and non-terminating nodes, along with their valuations. Other factorable polynomials are not as complete with their descriptions. 

Using these results, we can also look at specific trees and determine the polynomial of the lowest degree which could be represented by this tree, essentially the reverse process from what was detailed in this paper. The true generating polynomial would be a multiple of it where the other factors have a valuation of $0$.

Additionally, all the work done within this paper uses the field of rational numbers and $p$-adic valuations. However, other valuations and fields do exist. Whether or not anything remarkable can come from using these valuations with these trees remains to be seen.

\section*{Acknowledgments} 

This work was supported by the 2023 REU research grant at Ursinus College under NSF grant No. 1851948, for which I am grateful. I would also like to thank my advisor, Dr. Olena Kozhushkina at Ursinus College, for her guidance and mentoring, as well as Dr. Justin Trulen at Kentucky Wesleyan College for his advice and feedback. In addition, I am appreciative of my REU collaborators, MaeKayla Minton and Laura Vaughan, who encouraged and supported me through this process.



{\footnotesize  
\medskip
\medskip
\vspace*{1mm} 
 
\noindent {\it Dillon Snyder}\\  
University of Connecticut\\
341 Mansfield Road\\
Storrs, Connecticut\\
E-mail: {\tt dillon.snyder@uconn.edu}\\ \\  
  

}

\vspace*{1mm}\noindent\footnotesize{\date{ {\bf Received}: April 31, 2017\;\;\;{\bf Accepted}: June 31, 2017}}\\
\vspace*{1mm}\noindent\footnotesize{\date{  {\bf Communicated by Some Editor}}}


\newpage 

{\footnotesize
\begin{thebibliography}{00}
\bibitem{almo} L. Almodovar, A.N. Byrnes, J. Fink, X. Guan, A. Kesarwani, G. Lavigne, L.A. Medina, V.H. Moll, I. Nogues, S. Rajesekaran, E. Rowland, A. Yuan, A closed-form solution might be given by a tree. Valuations of quadratic polynomials, {\it Sci. Ser. A Math. Sci.}, {\bf 29} (2019), 11--28.

\bibitem{padic} W. Boultinghouse, E. Hammett, S. Hu, O. Kozhushkina, R. Snyder, J. Trulen, $p$-Adic valuations of quadratic sequences, {\it Preprint}, 1--19.

\bibitem{2adic} W. Boultinghouse, J. Long, O. Kozhushkina, J. Trulen, $2$-Adic valuations of quadratic sequences, {\it J. Integer Seq.}, {\bf 24} (2021), 1--22.

\bibitem{Gouvea} F. Gouv\^ea, {\it $p$-Adic Numbers: An Introduction}, Springer-Verlag, (1997).

\bibitem{guan} X. Guan, Methods in Symbolic Computation and $p$-Adic Valuations of Polynomials, available online at the URL: \url{https://ui.adsabs.harvard.edu/abs/2017PhDT........65G}

\bibitem{n2+7} O. Kozhushkina, M. Brucal-Hallare, J. Long, V.H. Moll, J.-C. Pedjeu, B. Thompson, J. Trulen, The valuation tree for $n^2+7$, {\it Sci. Ser. A Math. Sci.}, {\bf 30} (2020), 91--102.

\bibitem{medina} L.A. Medina, V.H. Moll, E. Rowland, Periodicity in the $p$-adic valuation of a polynomial, {\it J. Number Theory}, {\bf 180} (2017), 139--153.

\bibitem{nguyen} M. Nguyen, $2$-adic Valuations of Square Spiral Sequences, available online at the URL: \url{https://aquila.usm.edu/cgi/viewcontent.cgi?article=1774&context=honors_theses}



\end{thebibliography}}
\end{document}